\documentclass[12pt,leqno,twoside]{amsart}
\usepackage{amssymb,amsmath,amsthm,soul,color}
\usepackage{t1enc}
\usepackage[utf8]{inputenc}
\usepackage{a4,indentfirst,latexsym}
\usepackage{graphics}
\usepackage{mathrsfs}
\usepackage{cite,enumitem,graphicx}
\usepackage[colorlinks=true,urlcolor=blue,
citecolor=red,linkcolor=blue,linktocpage,pdfpagelabels,
bookmarksnumbered,bookmarksopen]{hyperref}
\usepackage[english]{babel}
\usepackage[left=2.61cm,right=2.61cm,top=2.72cm,bottom=2.72cm]{geometry}
\usepackage[colorinlistoftodos]{todonotes}
\usepackage{pifont}
\usepackage{mathrsfs}
\usepackage{cite,enumitem,graphicx}
\input xy
\xyoption{all}
\usepackage{color}

\newtheorem{Th}{Theorem}[section]

\newtheorem{Lem}[Th]{Lemma}

\newtheorem{Rem}[Th]{Remark}

\newenvironment{altproof}[1]
{\noindent
{\em Proof of {#1}}.}
{\nopagebreak\mbox{}\hfill $\Box$\par\addvspace{0.5cm}}

\makeatletter
    
    \newcommand{\Rmnum}[1]{\expandafter\@slowromancap\romannumeral #1@}

   \newcommand{\vp}{\varphi}

   \def\Z{\mathbb{Z}}

   \def\R{\mathbb{R}}

   \def\J{\mathcal{J}}

\newcommand{\cC}{{\mathcal C}}

\newcommand{\cI}{{\mathcal I}}
\newcommand{\cJ}{{\mathcal J}}

\newcommand{\cN}{{\mathcal N}}

\newcommand{\Ga}{\Gamma}

\newcommand{\weakto}{\rightharpoonup}

\numberwithin{equation}{section}

\begin{document}

\title{Nonlinear Schr\"odinger equations with sum of periodic and vanishing potentials and sign-changing nonlinearities}


\author[B. Bieganowski]{Bartosz Bieganowski}
\address[B. Bieganowski]{\newline\indent 
	Faculty of Mathematics and Computer Science,
	\newline\indent 
	Nicolaus Copernicus University,
	\newline\indent
	ul. Chopina 12/18, 87-100 Toru\'n, Poland}
\email{\href{mailto:bartoszb@mat.umk.pl}{bartoszb@mat.umk.pl}}

\author[J. Mederski]{Jaros\l aw Mederski}
\address[J. Mederski]{\newline\indent 
	Institute of Mathematics,
	\newline\indent
	Polish Academy of Sciences,
	\newline\indent 
	ul. \'Sniadeckich 8, 00-956
	Warszawa, Poland
	\newline\indent
	and
	\newline\indent
	Faculty of Mathematics and Computer Science,
	\newline\indent 
	Nicolaus Copernicus University,
	\newline\indent
	ul. Chopina 12/18, 87-100 Toru\'n, Poland}
\email{\href{mailto:jmederski@mat.umk.pl}{jmederski@impan.pl}}
\date{}
\maketitle

\pagestyle{myheadings} \markboth{\underline{B. Bieganowski, J. Mederski}}{
\underline{NLS equations with sum of periodic and vanishing potentials and sign-changing nonlinearities}}

\begin{abstract} We look for ground state solutions to the following nonlinear Schr\"odinger equation
$$-\Delta u + V(x)u = f(x,u)-\Gamma(x)|u|^{q-2}u\hbox{ on }\R^N,$$
where $V=V_{per}+V_{loc}\in L^{\infty}(\R^N)$ is the sum of a periodic potential $V_{per}$ and a localized potential $V_{loc}$, $\Gamma\in L^{\infty}(\R^N)$  is periodic and $\Gamma(x)\geq 0$ for a.e. $x\in\R^N$ and $2\leq q<2^*$. We assume that $\inf\sigma(-\Delta+V)>0$, where $\sigma(-\Delta+V)$ stands for the spectrum of $-\Delta +V$ and $f$ has the subcritical growth but higher than $\Gamma(x)|u|^{q-2}u$, however the nonlinearity $f(x,u)-\Gamma(x)|u|^{q-2}u$ may change sign. Although a Nehari-type monotonicity condition for the nonlinearity is not satisfied, we investigate the existence of ground state solutions being minimizers on the Nehari manifold.
\end{abstract}

\vspace{0.2cm}


\section*{Introduction}
\setcounter{section}{1}

We consider the following nonlinear Schr\"odinger equation
\begin{equation}
\label{eq}
-\Delta u + V(x)u = f(x,u)-\Gamma(x)|u|^{q-2}u\hbox{ on }\R^N,\; N\geq 1
\end{equation}
with $u\in H^1(\R^N)$, which appears in different areas of mathematical physics. In particular, if $V$, $f$ and $\Gamma$ are periodic (or close-to-periodic) in $x$, then there is a wide range of applications in photonic crystals admitting nonlinear effects \cite{Pankov,Kuchment}. In this case  \eqref{eq} describes the propagation of gap solitons which are special nontrivial solitary wave solutions $\Phi(x,t)=u(x)e^{-i\omega t}$ of the time-dependent Schr\"odinger equation of the form
\begin{equation*}\label{eq:SchrodTime}
i\frac{\partial \Phi}{\partial t}=
-\Delta \Phi+(V(x)+\omega)\Phi-g(x,|\Phi|),
\end{equation*}
where $g$ is responsible for nonlinear polarization in a photonic crystal, e.g 
in a self-focusing Kerr-like medium one has 
$f(x,u)=K(x) |u|^2 u$
for some periodic $K:\R^N\to\R$, $\inf K>0$ and $\Gamma\equiv0$; see \cite{Buryak,NonlinearPhotonicCrystals}. If $\Gamma(x)>0$ for $x\in\R^N$, then we deal with a mixture of self-focusing and defocusing optical materials. For instance, we consider $f(x,u)=K(x) |u|^{p-2}u$ with $2<q<p< 2^*$ and the nonlinear material may exhibit dual-power effect if $p=2q$.
In general $f$ satisfies the following conditions:

\begin{itemize}
\item[(F1)] $f:\R^N\times\R\to \R$ is measurable, $\Z^N$-periodic in $x\in\R^N$ and continuous in $u\in\R$ for a.e. $x\in\R^N$ and there are $c>0$ and $2\leq q<p< 2^*$ such that
$$|f(x,u)|\leq c(1+|u|^{p-1})\hbox{ for all }u \in\R,\; x\in\R^N.$$
\item[(F2)] $f(x,u)=o(|u|)$ uniformly in $x$ as $|u|\to 0$.
\item[(F3)] $F(x,u)/|u|^q\to\infty$ uniformly in $x$ as $|u|\to\infty$, where $F$ is the primitive of $f$ with respect to $u$.
\item[(F4)] $u\mapsto f(x,u)/|u|^{q-1}$ is strictly increasing on $(-\infty,0)$ and $(0,\infty)$.
\end{itemize} 

We impose on  $\Gamma$ the following condition: 
\begin{itemize}
\item[($\Gamma$)] $\Gamma\in L^{\infty}(\R^N)$ is periodic in $x\in\R^N$, $\Gamma(x)\geq 0$ for a.e. $x\in\R^N$.
\end{itemize} 
Observe that if $q=2$, then $\Gamma(x)|u|^{q-2}u=\Gamma(x)u$ may be included in the potential $V$ and we may assume that $\Gamma\equiv 0$ in this case.

In photonic crystals potential $V$ is periodic or close-to-periodic. Namely if the periodic structure has a linear defect, i.e. an additional structure breaking the periodicity, then the photonic crystal can guide light along the defect \cite{Doerfler,Kuchment}. In this case the potential has the following form 
\begin{equation}\label{PotentialForm}
V=V_{per}+V_{loc}, 
\end{equation} 
where $V_{per}$ is periodic in $x\in\R^N$ and $V_{loc}$ is a localized potential that vanishes at infinity; see condition (V).

Our goal is to find a {\em ground state} of the energy functional
$\J:H^1(\R^N)\to\R$ of class $\cC^1$
given by
$$\J(u)=\frac12\int_{\R^N}|\nabla u|^2+V(x)|u|^2\,dx- \int_{\R^N} \Big(F(x,u)-\frac1q\Gamma(x)|u|^q\Big)\, dx.$$
A ground state stands for a critical point  being a minimizer of $\J$ on the Nehari manifold 
$$\cN:=\{u\in H^1(\R^N)\setminus\{0\}:\; \J'(u)(u)=0\}.$$
Obviously $\cN$ contains all nontrivial critical points, hence a ground state is the least energy solution.

Recall that in the absence of the localized potential $V_{loc}=0$, the spectrum $\sigma(-\Delta+V)$ of $-\Delta+V=-\Delta+V_{per}$ is purely continuous, bounded from below and consists of closed disjoint intervals \cite{ReedSimon}. In this case ground states of \eqref{eq} with $\Gamma\equiv0$ has been of particular interests of many authors; see e.g.
\cite{Rabinowitz:1992,LiWangZeng,Pankov,SzulkinWeth,Liu,MederskiNLS2014}
and references therein. In view of a recent result of Szulkin and Weth \cite{SzulkinWeth}, problem (\ref{eq}) under assumptions (F1)-(F4) with $q=2$, $2<p<2^*$, $\Gamma\equiv 0$ and $V=V_{per}\in \cC(\R^N,\R)$, $0\notin\sigma(-\Delta+V)$ admits a minimizer on the Nehari-Pankov manifold \cite{Pankov} which is, in general, contained in $\cN$ but coincides with $\cN$ provided that $0<\inf\sigma(-\Delta+V)$. 
Note that if $\Gamma\neq 0$ and $q>2$, then the nonlinear part of the energy functional 
$$\cI(u):=\int_{\R^N} \Big(F(x,u)-\frac1q\Gamma(x)|u|^q\Big)\, dx$$
is sign-changing, moreover $u\mapsto \big(f(x,u)-\Gamma(x)|u|^{q-2}u\big)/|u|$ is no longer increasing on $(-\infty,0)$ and $(0,\infty)$, so the results of \cite{SzulkinWeth} do not apply in our case. Moreover under our assumptions, $\cN$ is not $\cC^1$-manifold, so that the classical minimization on the Nehari manifold does not work. 
We intend to adopt the techniques of \cite{SzulkinWeth} based on the observation that $\cN$ is a topological manifold homeomorphic with the unit sphere in $H^1(\R^N)$, where a minimizing sequence can be found. Our approach is presented in the abstract setting in Section \ref{sect:VariationalSetting} and 
we develop a critical point theory which extends the abstract result from \cite{BartschMederski1} for positive definite functionals 
and enables us to deal with sign-changing nonlinearities; see Theorem \ref{ThSetting}. Note that abstract results concerning Nehari techniques have been also obtained by Szulkin and Weth in \cite{SzulkinWethHandbook} for positive nonlinear part $\cI$ or completely continuous $\cI'$  as well as by Figueiredo and Quoirin in \cite{Figueiredo} for weakly lower semicontinuous $u\mapsto\J'(u)u$. Observe that in our problem $\cI$ and $\J$ do not satisfy these conditions anymore. 

Sign-changing nonlinearities of the form $a(x)g(u)$, where $a(x)$ changes sign (called indefinite nonlinearities), have been studied for instance in \cite{CostaTehrani, LiuYang, CostaTehraniCalVPDE}. In our case, however, the sign of the nonlinear term depends on $u$, and the results devoted to indefinite nonlinearities do not apply.

In this paper we consider a potential of the form \eqref{PotentialForm} and
we work with the following assumption:
\begin{itemize}
\item[($V$)] $V_{per}\in L^{\infty}(\R^N)$ is $\Z^N$-periodic, $V_{loc}\in L^{\infty}(\R^N)$ and $V_{loc}(x)\to 0$ as $|x|\to\infty$.
\end{itemize}
Observe that for the sufficiently decaying perturbation $V_{loc}$, the multiplication by
$V_{loc}$ is a compact perturbation of $-\Delta+V_{per}$; see \cite{ReedSimon,PankovNotes}.
Therefore the essential spectrum 
$$\sigma_{ess}(-\Delta+V)=\sigma_{ess}(-\Delta+V_{per})=\sigma(-\Delta+V_{per}).$$  
However the whole spectrum $\sigma(-\Delta+V)$ is no longer purely continuous and may contain eigenvalues below the essential part $\sigma_{ess}(-\Delta+V)=\sigma(-\Delta+V_{per})$. 

Now we state our main results.
\begin{Th}\label{ThMain1}
Suppose that (V), ($\Gamma$) (F1)-(F4) are satisfied and  $\inf\sigma(-\Delta+V)>0$. If $V_{loc}(x) < 0$ for a.e. $x\in\R^N$ or $V_{loc} \equiv 0$, then \eqref{eq} has a ground state $u\in H^1(\R^N)$, i.e. $u$ is a critical point of $\J$ such that $\J(u)=\inf_{\cN}\J$. Moreover $u$ is continuous and there exist $\alpha, C>0$ such that
$$|u(x)|\leq C \exp(-\alpha |x|)\hbox{ for any }x\in\R^N.$$
\end{Th}

\begin{Th}\label{ThNonexistence}
Suppose that (V), ($\Gamma$) (F1)-(F4) are satisfied and  $\inf\sigma(-\Delta+V_{per})>0$. If $V_{loc}(x) > 0$ for a.e. $x \in \R^N$ 
then \eqref{eq} has no ground states.
\end{Th}

In order to find ground states we develope a critical point theory in Section \ref{sect:VariationalSetting}  and prove Theorems \ref{ThMain1} and \ref{ThNonexistence} in Section \ref{sect:Appliaction}. Due to a localized potential $V_{loc}$ we need to investigate decompositions of Palais-Smale sequences and we present this result in a self-contained Section \ref{sect:Decomposition}.

\section{Abstract setting}
\label{sect:VariationalSetting}

Let $E$ be a Hilbert space with the norm $\| \cdot \|$. We consider a functional $\mathcal{J} : E \rightarrow \R$ of the following form
$$
\cJ (u) = \frac{1}{2} \|u\|^2 - \cI(u),
$$
where $\cI : E \rightarrow \R$ is of $\cC^1$ class. In this case the Nehari manifold is given by
\begin{eqnarray*}
\cN &:=& \{ u \in E \setminus \{0\} \ : \ \cJ'(u)(u) = 0 \} \\ &=& \{ u \in E \setminus \{0\} \ : \ \|u\|^2 = \cI'(u)(u) \}.
\end{eqnarray*}

Now we formulate the main result of this section.

\begin{Th}\label{ThSetting}
Suppose that the following conditions hold:
\begin{itemize}
	\item[(J1)] there is $r > 0$ such that $a:= \inf_{\|u\|=r} \cJ(u) > \J(0)=0$; 
	\item[(J2)] there is $q\geq 2$ such that $\cI(t_nu_n)/t_n^q \to \infty$ for any $t_n\to\infty$ and $u_n\to u\neq 0$ as $n\to\infty$;
	\item[(J3)]  for $t \in (0,\infty)\setminus \{1\}$ and $u\in\cN$ 
	$$\frac{t^2-1}{2} \cI'(u)(u)-\cI(tu)+\cI(u) < 0;$$ 
	\item[(J4)] $\cJ$ is coercive on $\cN$.
\end{itemize}
Then $\inf_{\cN}\J>0$ and there exists a bounded minimizing sequence for $\cJ$ on $\cN$, i.e. there is a  sequence $(u_n) \subset \cN$ such that $\cJ(u_n) \to \inf_{\cN} \cJ$ and $\cJ'(u_n) \to 0$.
\end{Th}
Observe that condition (J2) implies that for any $u\neq 0$ there is $t>0$ such that $\J(tu)<0$, hence taking into account also (J1) we easily check that $\J$ has the classical mountain pass geometry \cite{AmbrosettiRabinowitz,Willem} and we are able to find a Palais-Smale sequence. However we do not know whether it is a bounded sequence and contained in $\cN$. In order to get the boundedness we assume the coercivity in (J4), which is, in applications, a weaker requirement than the classical Ambrosetti-Rabinowitz condition; see e.g. \cite{SzulkinWeth}. 

\begin{Rem}\label{remark} a) In order to get (J3) it is sufficient to check 
\begin{eqnarray}
(1-t) \left( t\cI'(u)(u) - \cI'(tu)(u)\right) > 0&&\hbox {for any }t \in (0, \infty) \setminus \{1\}  \hbox{ and }\\&&u\hbox{ such that }\cI'(u)(u)>0.\nonumber
\end{eqnarray}
Indeed,
let us consider $t \in (0,\infty)\setminus \{1\}$, $u\in\cN$ and
\begin{equation}\label{eq:defphi}
\varphi(t) = \frac{t^2-1}{2} \cI'(u)(u)-\cI(tu)+\cI(u).
\end{equation}
Then $\cI'(u)(u)=\|u\|^2>0$, $\varphi(1) = 0$, $\varphi'(t) = t \cI'(u)(u) - \cI'(tu)(u)>0$ for $t<1$ and $\varphi'(t)<0$ for $t>1$.
Therefore $\varphi(t) < \varphi(1) = 0$.\\
b) Observe that (J3) is equivalent to the following condition: $u\in \cN$ is the unique maximum point of $ (0,+\infty)\ni t\mapsto \J(tu)\in\R$. Indeed, note that for $u\in\cN$
\begin{equation}\label{eq:ineqJ}
\J(tu)=\J(u)+\big(\J(tu)-\J(u)-\frac{t^2-1}{2}\J'(u)(u)\big)=\J(u)+\vp(t)<\J(u)
\end{equation}
if and only if $\vp(t)<0$.
\end{Rem}

\begin{altproof}{Theorem \ref{ThSetting}}
For a given $u \neq 0$ we consider a map $\vp:[0,+\infty)\to\R$ defined by $\varphi(t) = \cJ(tu) - \cJ(u)$ for $t\in [0,+\infty)$. Note that from \eqref{eq:ineqJ}, $\vp(t)$ is given by \eqref{eq:defphi} provided that $u\in\cN$. In view of
(J1)-(J2) we obtain
$\vp(0)=-\J(u)<\vp(\frac{r}{\|u\|})$ and
 $\vp(t)\to-\infty$ as $t\to\infty$. Therefore there is a maximum point $t(u)>0$ of $\vp$ which is a critical point of $\varphi$, i.e. $\cJ'(t(u) u)(u) = 0$ and 
 $t(u) u \in \cN$.  In view of Remark \ref{remark} b) we infer that for any $u\neq 0$ there is an unique critical point $t(u)>0$ of $\vp$, i.e. $t(u) u \in \cN$.
Let $\hat{m} : E \setminus \{0\} \rightarrow \cN$ be a map given by $\hat{m}(u) = t(u) u$ for $u\neq 0$. We are going to show that $\hat{m}$ is continuous. Take $u_n\to u_0\neq 0$ and denote $t_n=t(u_n)$ for $n\geq 0$, so that $\hat{m}(u_n)=t_nu_n$. Observe that if $t_n\to\infty$ then by \eqref{eq:ineqJ} and (J2) 
$$o(1)=\J(u_n)/t_n^q\leq \J(\hat{m}(u_n))/t_n^q=\frac12\|u_n\|^2t_n^{2-q}-\cI(t_nu_n)/t_n^q\to-\infty$$
and we get a contradiction. Therefore we may assume that $t_n\to t_0\geq 0$. Again by \eqref{eq:ineqJ}
$$\J(t_{u_0}u_0)\geq \J(t_0u_0)=\lim_{n\to\infty}\J(t_nu_n)\geq \lim_{n\to\infty} \J(t_{u_0}u_n)=\J(t(u_0)u_0)$$
we get $t_0=t(u_0)$, which completes the proof of continuity of $\hat{m}$.
Thus $m = \hat{m}|_{S^1}$, where $S^1$ is the unit sphere in $E$, is a homeomorphism. The inverse function is given by $m^{-1}(v) = v/\|v\|$ for $v\in\cN$. Therefore
$$
c:=\inf_{u\in S^1}(\cJ \circ m)(u) = \inf_{u \in \cN} \cJ(u) \geq \inf_{u \in \cN} \cJ \left( \frac{r}{\|u\|}u \right) \geq a > 0.
$$
Now, arguing as in \cite{SzulkinWeth}[Proposition 2.9], we show  that $\cJ \circ \hat{m}$ is of $\cC^1$-class. Moreover for $w,z \in E$, $w \neq 0$ we have
$$
(\cJ \circ \hat{m})'(w)z = \frac{\|\hat{m}(w)\|}{\|w\|} \cJ '(\hat{m}(w)) z.
$$
Since $\cJ \circ m : S^1 \rightarrow \R$ is of $\cC^1$-class, then in view of the Ekeland variational principle we find a minimizing sequence $(v_n)\subset S^1$ for $\cJ \circ m$ such that $(\cJ \circ m)'(v_n)\to 0$. Then take  $u_n=m(v_n)\in\cN$ and observe that $\J'(u_n)(v_n)=0$ and
$$
(\cJ \circ m) ' (v_n) (z) = \|u_n\| \cJ ' (u_n)(z) = \|u_n\| \cJ ' (u_n)(z+tv_n).
$$
for any $z\in T_{v_n}S^1$ and $t\in\R$, where $T_{v_n}S^1$ stands for the tangent space $S^1$ at $v_n$. 
Therefore
$$
\|(\J\circ m)'(v_n)\|=\sup_{z \in T_{v_n} S^1, \ \|z\|=1} (\cJ \circ m) ' (v_n) (z) = \|u_n\| \|\cJ ' (u_n) \|.
$$
Since $u_n \in \cN$, we have $\|u_n\| \geq \eta$ for some $\eta>0$. The coercivity of $\cJ$ implies that $\sup_n \|u_n\| < \infty$. Hence
$(u_n)$ is a bounded minimizing  sequence for $\cJ$ on $\cN$ such that $\J'(u_n)\to 0$.
\end{altproof}

\section{Application to \eqref{eq}}
\label{sect:Appliaction}

\begin{Lem}
Suppose that (V), ($\Ga$), (F1)-(F4) are satisfied and $\inf\sigma(-\Delta+V)>0$. Then (J1)-(J4) hold.
\end{Lem}
In $H^1 (\R^N)$ we consider the following norm
$$
\| u \|^2 = \int_{\R^N} | \nabla u |^2 \, dx + \int_{\R^N} V(x) u^2 \, dx.
$$
Since $\inf\sigma(-\Delta+V)>0$, the norm $\| \cdot \|$ is equivalent to the classic one on $H^1 (\R^N)$.
\begin{proof}
\begin{itemize}
	\item[(J1)] Fix $\varepsilon > 0$. Observe that (F1) and (F2) implies that $F(x,u) \leq \varepsilon |u|^2 + C_\varepsilon |u|^p$ for some $C_\varepsilon$. Therefore 
	$$
	\int_{\R^N} F(x,u) \, dx - \int_{\R^N} \frac{1}{q} \Ga (x) |u|^q \, dx \leq \int_{\R^N} F(x,u) \, dx \leq C(\varepsilon \|u\|^2 + C_\varepsilon \|u\|^p),
	$$
for some constant $C > 0$ provided by the Sobolev embedding theorem. Thus there is $r > 0$ such that 
$$
\int_{\R^N} F(x,u) \, dx - \int_{\R^N} \frac{1}{q} \Ga (x) |u|^q \, dx \leq \frac{1}{4} \|u\|^2
$$
for $\|u\| \leq r$.
Therefore
$$
\cJ (u) \geq \frac{1}{4} \|u\|^2 = \frac{1}{4} r^2 > 0
$$
for $\|u\| = r$.
	\item[(J2)] By (F3) and Fatou's lemma we get
	$$
	\cI (t_n u_n) / t_n^q = \int_{\R^N} \frac{F(x, t_n u_n)}{t_n^q} \, dx - \frac{1}{q} \int_{\R^N} \Ga (x) |u_n|^q \, dx \to \infty.
	$$
	\item[(J3)] Fix $u$ such that $\cI'(u)(u)>0$, i.e. 
	$$
\int_{\R^N} f(x,u)u \, dx > \int_{\R^N} \Ga(x) |u|^q \, dx.
$$
Observe that
\begin{eqnarray*}
t \cI'(u)(u) - \cI ' (tu)(u) &=& \int_{\R^N} f(x,u)tu - \Ga(x) t |u|^q - f(x,tu)u + \Ga(x) t^{q-1} |u|^q \, dx\\
&=& \int_{\R^N} f(x,u)tu - f(x,tu)u \, dx + (t^{q-1}-t) \int_{\R^N} \Ga (x) |u|^q \, dx.
\end{eqnarray*}
Then for $t < 1$ we have
\begin{eqnarray*}
& & \int_{\R^N} f(x,u)tu - f(x,tu)u \, dx + (t^{q-1}-t) \int_{\R^N} \Ga (x) |u|^q \, dx > \\
&&> \int_{\R^N} t^{q-1} f(x,u)u - f(x,tu)u \, dx = \\
&&= t^{q-1} \int_{\R^N} f(x,u)u - \frac{f(x,tu)u}{t^{q-1}} \, dx > 0,
\end{eqnarray*}
by (F4). In a similar way we see that $t \cI'(u)(u) - \cI ' (tu)(u) < 0$ for $t>1$. In view of Remark \ref{remark}, the condition (J3) holds.

\item[(J4)] Let $(u_n) \subset \cN$ be a sequence such that $\|u_n\| \to \infty$ as $n\to\infty$. If $q=2$, then  $\Gamma\equiv 0$. Suppose that 
$v_n:=u_n(\cdot+y_n)$ and if $v_n\to0$ in $L^p(\R^N)$, then by \eqref{eq:ineqJ} one has
$$\J(u_n)\geq\J(sv_n)=\frac{s^2}{2}+o(1)$$
for any $s>0$. Hence $\J(u_n)\to\infty$. If $v_n$ is bounded away from $0$ in $L^p(\R^N)$, then applying Lions' lemma and Fatou's lemma we show that
$$\J(u_n)/\|u_n\|^2=\frac12-\int_{\R^N}\frac{F(x,u_n(x+y_n))}{u_n^2(x+y_n)}v_n^2(x+y_n)\,dx\to-\infty$$
for some sequence $(y_n)\subset\R^N$ such that $v_n(\cdot+y_n)\weakto v\neq 0$ for some $v\in H^1(\R^N)$ and $v_n(x+y_n)\to v(x)$  for a.e. $x\in\R^N$. Thus we get a contradiction since $\J(u_n)/\|u_n\|^2\geq 0$. If $q>2$ then (F3) implies that
$$
f(x,u)u = q\int_0^u \frac{f(x,u)}{u^{q-1}}s^{q-1}\,ds \geq q\int_0^u \frac{f(x,s)}{s^{q-1}}s^{q-1}\,ds=q F(x,u)
$$
for $u\geq 0$ and similarly $f(x,u)u\geq q F(x,u)$ for $u<0$.
Therefore
\begin{eqnarray*}
\cJ(u_n) &=& \frac{1}{2} \|u_n\|^2 - \int_{\R^N} F(x, u_n) \, dx + \frac{1}{q} \int_{\R^N} \Ga (x) |u_n|^q \, dx = \\
&=& \left( \frac{1}{2} - \frac{1}{q} \right) \|u_n\|^2 + \int_{\R^N} \frac{1}{q} f(x, u_n) u_n - F(x, u_n) \, dx \geq \\
&\geq& \left( \frac{1}{2} - \frac{1}{q} \right) \|u_n\|^2 \to \infty
\end{eqnarray*}
\end{itemize}
as $n\to\infty$.
\end{proof}


\begin{altproof}{Theorem \ref{ThMain1}}
In view of  Theorem \ref{ThSetting} we find a bounded minimizing sequence $(u_n) \subset \cN$ for $\cJ$ on $\cN$. Let $c:=\inf_{\cN}\J$, 
$$
\cJ_{per} (u) := \frac{1}{2} \int_{\R^N} |\nabla u|^2 + V_{per}(x) u^2 \, dx - \int_{\R^N} F(x,u) \, dx + \frac{1}{q} \int_{\R^N} \Ga(x) |u|^q \, dx
$$
and
$$
c_{per} := \inf \{ \cJ_{per} (u) \ : \ u \in H^1 (\R^N) \setminus \{0\}, \ \cJ_{per}' (u) = 0 \}.
$$
If $V_{loc} = 0$, we have $\cJ = \cJ_{per}$. Therefore, in view of Theorem \ref{ThDecomposition}, passing to a subsequence,
there is  an integer $\ell \geq 0$ and sequences $(y_n^k) \subset \mathbb{Z}^N$, $w^k \in H^1 (\R^N)$, $k = 1, \ldots, \ell$ such that $u_n \rightharpoonup u_0$, $\cJ' (u_0) = 0$,
$w^k \neq 0$, $\cJ_{per}'(w^k) = 0$ for each $1 \leq k \leq \ell$ and
$$
\cJ_{per} (u_n) \to \cJ_{per} (u_0) + \sum_{k=1}^\ell \cJ_{per} (w^k).
$$
Thus
$$
c = c_{per} = \cJ_{per} (u_0) + \sum_{k=1}^\ell \cJ_{per} (w^k) \geq \cJ_{per} (u_0) + \ell c_{per}.
$$
If $u_0 \neq 0$, then $\ell=0$ and $\cJ_{per} (u_0)=c_{per}$ and $u_0$ is a ground state. If $u_0 = 0$, we get $c_{per} \geq \ell c_{per}$. Since $c_{per} > 0$ then $\ell = 1$ and $w^1$ is a ground state. \\
Suppose that $V_{loc} < 0$ for a.e. $x \in \R^N$. Again, using Theorem \ref{ThDecomposition}, there is $u_{per} \neq 0$ - a critical point of $\cJ_{per}$ such that $\cJ_{per} (u_{per}) = c_{per}$. Let $t > 0$ be such that $tu_{per} \in \cN$. Since $V(x) < V_{per} (x)$ we get
$$
c_{per} = \cJ_{per} (u_{per}) \geq \cJ_{per} (tu_{per}) > \cJ (tu_{per}) \geq c > 0.
$$
From Theorem \ref{ThDecomposition} there is $u_0$ such that $\cJ ' (u_0) = 0$. Moreover
$$
\cJ (u_n) \to \cJ (u_0) + \sum_{k=1}^\ell \cJ_{per} (w^k).
$$
Thus
$$
c = \cJ (u_0) + \sum_{k=1}^\ell \cJ_{per} (w^k) \geq \cJ (u_0) + \ell c_{per}.
$$
Since $c_{per} > c$, we get that $\ell = 0$ and $u_0$ is a ground state. From Theorem 2 in \cite{PankovDecay} we know that there exist $\alpha, C > 0$ such that
$$
|u(x)| \leq C \exp (-\alpha |x|)
$$
for $x \in \R^N$.
\end{altproof}





\subsection{Nonexistence result}
\label{sect:Nonexistence}

We will prove that if $V_{loc} (x) > 0$ for a.e. $x \in \R^N$, then ground states  do not exist. Observe that the condition $\inf(-\Delta+V_{per})>0$ implies that the norm of given by $\big(\int_{\R^N} |\nabla u |^2 + V_{per}(x) u^2 \, dx\big)^{1/2}$ for $u\in H^1(\R^N)$, is equivalent to  $\|\cdot\|$.

\begin{altproof}{Theorem \ref{ThNonexistence}}
Let us denote
\begin{eqnarray*}
\cN_{per} &=& \left\{ u \in H^1 (\R^N) \setminus \{0\} \ : \J'_{per}(u)(u)=0 \right\}
\end{eqnarray*}
where
$$
\cJ_{per} (u) = \cJ (u) - \int_{\R^N} V_{loc}(x) u^2 \, dx.
$$
Suppose, by contradiction, that there is a ground state $u_0 \in \cN$ of $\J$. Let $t_{per} > 0$ be such that $t_{per} u_0 \in \cN_{per}$. Since $V_{loc}(x) > 0$ for a.e. $x \in \R^N$, we have that
$$
\int_{\R^N} V_{loc}(x) u_0^2 \, dx > 0.
$$
Therefore
$$
c_{per} := \inf_{\cN_{per}} \cJ_{per} \leq \cJ_{per} (t_{per} u_0) < \cJ(t_{per} u_0) \leq \cJ(u_0) = c.
$$
On the other hand, let $u \in \cN_{per}$ and for $y \in \mathbb{Z}^N$ let us denote $u_y = u(\cdot - y)$. For each $y \in \mathbb{Z}^N$ let $t_y$ be a number such that $t_y u_y \in \cN$. Now observe that
$$
\cJ_{per} (u) = \cJ_{per}(u_y) \geq \cJ_{per}(t_y u_y) = \cJ(t_y u_y) - \int_{\R^N} V_{loc}(x) (t_y u_y)^2 \, dx \geq  -c \int_{\R^N} V_{loc}(x) (t_y u_y)^2 \, dx.
$$
We are going to show that 
$$
\int_{\R^N} V_{loc}(x) (t_y u_y)^2 \, dx \to 0.
$$
Indeed, note that
$$
\int_{\R^N} V_{loc}(x) (t_y u_y)^2 \, dx = t_y^2 \int_{\R^N} V_{loc}(x + y) u^2 \, dx.
$$
In view of (V) we easily get
$$
\int_{\R^N} V_{loc}(x + y) u^2 \, dx \to 0
$$
as $|y|\to\infty$.
Since $\cJ_{per}$ is coercive on $\mathcal{N}_{per}$, then
$$
\cJ_{per} (t_y u_y) = \cJ_{per} (t_y u) \leq c_{per}
$$
implies that $(t_y)$ is bounded. Therefore
$$
\cJ_{per} (u) \geq c + \int_{\R^N} V_{loc}(x) (t_y u_y)^2 \, dx \to c.
$$
Taking infimum over all $u \in \cN_{per}$ we have a contradiction $c_{per} \geq c$.
\end{altproof}

\section{Decomposition of bounded Palais-Smale sequences}
\label{sect:Decomposition}

In this section we provide a decomposition result of bounded Palais-Smale sequences in the spirit of \cite{Lions84,CotiZelati} which is a key step in proof of Theorem \ref{ThMain1} and generalizes \cite{JeanjeanTanakaIndiana2005}[Theorem 5.1] by Jeanjean and Tanaka. A similar result has been obtained for strongly indefinite functionals involving a Hardy potential and a more restrictive class of nonlinearities in \cite{GuoMederski}.

Let us consider the functional $\cJ : H^1 (\R^N) \rightarrow \R$ of the form
$$
\cJ(u) = \frac{1}{2} \int_{\R^N} |\nabla u|^2 + V(x) u^2 \, dx - \int_{\R^N} G(x,u) \, dx.
$$
We assume that (V) holds, $\inf \sigma (-\Delta + V(x)) > 0$ and $G(x,u) = \int_0^u g(x,s) \, ds$, where $g:\R^N\times \R\to\R$ satisfies the following conditions:
\begin{itemize}
\item[($G1$)] $g(\cdot, u)$ is measurable and $\mathbb{Z}^N$-periodic in $x$, $g(x,\cdot)$ is continuous in $u$ for a.e. $x \in \R^N$;
\item[($G2$)] $g(x,u) = o(|u|)$ as $|u| \to 0$ uniformly in $x$;
\item[($G3$)] there exists $2 < r < 2^*$ such that $\lim_{|u| \to \infty} g(x,u)/|u|^{r-1} = 0$ uniformly in $x$;
\item[($G4$)] for each $a<b$ there is a constant $c>0$ such that  $|g(x,u)|\leq c$ for a.e. $x\in\R^N$ and $a\leq u\leq b$.
\end{itemize} 

Let $\cJ_{per}:H^1 (\R^N)\to\R$ be given by the following formula
$$
\cJ_{per}(u) = \frac{1}{2} \int_{\R^N} |\nabla u|^2 + V_{per}(x) u^2 \, dx - \int_{\R^N} G(x,u) \, dx.
$$
Then we have the following result.
\begin{Th}\label{ThDecomposition}
Suppose that ($G1$)-($G4$), ($V1$), ($V2$) hold. Let $(u_n)$ be a bounded Palais-Smale sequence for $\cJ$. Then passing to a subsequence of $(u_n)$, there is an integer $\ell > 0$ and sequences $(y_n^k) \subset \mathbb{Z}^N$, $w^k \in H^1 (\R^N)$, $k = 1, \ldots, \ell$ such that:
\begin{itemize}
\item[$(a)$] $u_n \rightharpoonup u_0$ and $\cJ' (u_0) = 0$;
\item[$(b)$] $|y_n^k| \to \infty$ and $|y_n^k - y_n^{k'}| \to \infty$ for $k \neq k'$;
\item[$(c)$] $w^k \neq 0$ and $\cJ_{per}'(w^k) = 0$ for each $1 \leq k \leq \ell$;
\item[$(d)$] $u_n - u_0 - \sum_{k=1}^\ell w^k (\cdot - y_n^k) \to 0$ in $H^1 (\R^N)$ as $n \to \infty$;
\item[$(e)$] $\cJ (u_n) \to \cJ(u_0) + \sum_{k=1}^\ell \cJ_{per} (w^k)$.
\end{itemize}
\end{Th}

Observe that ($G2$) -- ($G4$) imply that for each $\varepsilon > 0$ there is $C_\varepsilon > 0$ such that 
\begin{equation}\label{eq:NonIneq}
|g(x,u)| \leq \varepsilon |u| + C_\varepsilon |u|^{r-1}
\end{equation}
for a.e. $x\in\R^N$ and for any $u\in \R$. In fact, we use only \eqref{eq:NonIneq} in proof of Theorem \ref{ThDecomposition} instead of (G2) -- (G4). (G1) -- (G4) can be easy verified in applications, e.g.  in our problem \eqref{eq} with $g(x,u)=f(x,u)-\Gamma(x)|u|^{q-2}u$. Similarly as in \cite{JeanjeanTanakaIndiana2005}, where an autonomous nonlinearity has been considered, our argument is based on Lions' lemma \cite{Lions84,Willem} and for the reader's conveniece we provide all details. In particular, Step 6 requires some technical arguments, since our nonlinearity is non-autonomous and the periodic part $V_{per}$ of potential $V$ is not constant.

Note that, in view of the Sobolev embeddings  and \eqref{eq:NonIneq} there exists $\rho > 0$ such that each nontrivial, critical point $v$ of $\cJ$ satisfies $\|v\| \geq \rho$, and now we are ready to prove the decomposition result.

\begin{altproof}{Theorem \ref{ThDecomposition}}  \\
\textbf{Step 1:} \textit{We may find a subsequence of $(u_n)$ such that $u_n \rightharpoonup u_0$, where $u_0 \in H^1 (\R^N)$ is a critical point of $\cJ$.} \\
Indeed, $(u_n)$ is bounded and therefore there is $u_0 \in H^1 (\R^N)$ such that $u_n \rightharpoonup u_0$ up to a subsequence; moreover, up to a subsequence, we may assume that $u_n (x) \to u(x)$ for a.e. $x \in \R^N$. Let $\varphi \in \cC_0^\infty (\R^N)$. Then
\begin{eqnarray*}
\cJ'(u_n)(\varphi) - \cJ'(u_0)\varphi &=& \int_{\R^N} \nabla (u_n - u_0) \nabla \varphi \, dx + \int_{\R^N} V(x) (u_n - u_0) \varphi \, dx \\
&-& \int_{\R^N} \left( g(x,u_n)-g(x,u_0) \right) \varphi \, dx \\
&=& \int_{\mathrm{supp}\, \varphi} \nabla (u_n - u_0) \nabla \varphi \, dx + \int_{\mathrm{supp}\, \varphi} V(x) (u_n - u_0) \varphi \, dx \\
&-& \int_{\mathrm{supp}\, \varphi} \left( g(x,u_n)-g(x,u_0) \right) \varphi \, dx.
\end{eqnarray*}
In view of the weak convergence $u_n \rightharpoonup u_0$ we see that
$$
\int_{\mathrm{supp}\, \varphi} \nabla (u_n - u_0) \nabla \varphi \, dx \to 0.
$$
By the Vitali convergence theorem we have
$$
\int_{\mathrm{supp}\, \varphi} V(x) (u_n - u_0)  \varphi \, dx \to 0.
$$
Next, the H\"older inequality and \eqref{eq:NonIneq} imply
$$
\int_{\mathrm{supp}\, \varphi} | g(x,u_n) \varphi| \, dx \leq C ( |u_n|_2 |\varphi \chi_E|_2 + |u_n|_r^{r-1} |\varphi \chi_E|_r )
$$
for measurable set $E \subset \mathrm{supp}\, \varphi$ and therefore
$$
\int_{\mathrm{supp}\, \varphi} \left( g(x,u_n)-g(x,u_0) \right) \varphi \, dx \to 0.
$$
Hence $\cJ'(u_n)\varphi \to \cJ'(u_0)\varphi$. Since $(u_n)$ is a Palais-Smale sequence, then $\cJ'(u_0) \varphi = 0$ for any $\varphi\in\cC_0^{\infty}(\R^N)$. \\
\textbf{Step 2:} \textit{Let $v_n^1 = u_n - u_0$. Suppose that}
\begin{equation}\label{eq:step2}
\sup_{z \in \R^N} \int_{B(z,1)} |v_n^1|^2 \, dx \to 0.
\end{equation}
\textit{Then $u_n \to u_0$ and (a)--(e) hold for $\ell = 0$.} \\
Note that
\begin{eqnarray*}
\cJ'(u_n)(v_n^1) &=& \int_{\R^N} \nabla v_n^1 \nabla v_n^1 \, dx + \int_{\R^N} \nabla u_0 \nabla v_n^1 \, dx + \int_{\R^N} V(x) |v_n^1|^2 \, dx \\
&+& \int_{\R^N} V(x) u_0 v_n^1 \, dx - \int_{\R^N} g(x,u_n) v_n^1 \, dx.
\end{eqnarray*}
Therefore
$$
\|v_n^1 \|^2 = \cJ'(u_n)(v_n^1) - \int_{\R^N} \nabla u_0 \nabla v_n^1 \, dx - \int_{\R^N} V(x) u_0 v_n^1 \, dx + \int_{\R^N} g(x,u_n) v_n^1 \, dx.
$$
Since $\cJ'(u_0)(v_n^1) = 0$ we get
$$
\|v_n^1 \|^2 = \cJ'(u_n)(v_n^1) + \int_{\R^N} \left( g(x,u_n) - g(x,u_0) \right) v_n^1 \, dx.
$$
From Lions' lemma (see \cite{Willem,Lions84}), we have $v_n^1 \to 0$ in $L^s (\R^N)$ for each $2<s<2^*$. Since $(v_n^1)$ is bounded, we observe that
$$
\| \cJ'(u_n)(v_n^1) \| \leq \| \cJ'(u_n) \| \| v_n^1 \| \to 0.
$$
Moreover, in view of a H\"older inequality and \eqref{eq:NonIneq}
$$
\left| \int_{\R^N} g(x,u_n) v_n^1 \, dx \right| \leq \varepsilon |u_n|_2 |v_n^1|_2 + C_\varepsilon |u_n|_{r}^{r-1} |v_n^1|_r.
$$
Therefore $\int_{\R^N} g(x,u_n) v_n^1 \, dx \to 0$ and in a very similar way $\int_{\R^N} g(x,u_0) v_n^1 \, dx \to 0$. Thus $v_n^1 \to 0$ and this completes the proof of Step 2. \\
\textbf{Step 3:} \textit{Suppose that there is a sequence $(z_n) \subset \mathbb{Z}^N$ such that
$$
\liminf_{n\to\infty} \int_{B(z_n, 1+\sqrt{N})} |v_n^1|^2 \, dx > 0.
$$
Then there is $w \in H^1 (\R^N)$ such that (up to a subsequence):}
$$
(i) \ |z_n| \to \infty, \quad (ii) \ u_n(\cdot + z_n) \rightharpoonup w \neq 0, \quad (iii) \ \cJ_{per}' (w) = 0.
$$
Conditions $(i)$ and $(ii)$ are standard, so let us concentrate on $(iii)$. Let $v_n = u_n(\cdot + z_n)$, then similarly as in Step 1
$$
\cJ_{per}'(v_n)(\varphi) - \cJ_{per}'(w)(\varphi) \to 0
$$
for each $\varphi \in \cC_0^\infty (\R^N)$. We are going to show that $\cJ_{per}'(v_n)(\varphi) \to 0$. Observe that
\begin{eqnarray*}
\cJ'(u_n) (\varphi(\cdot - z_n)) &=& \int_{\R^N} \nabla u_n(\cdot + z_n) \nabla\varphi \, dx + \int_{\R^N} V(x+z_n) u_n (\cdot + z_n) \, dx \\
&-& \int_{\R^N} g(x, u_n(\cdot +z_n)) \varphi \, dx.
\end{eqnarray*}
Moreover, there is $D > 0$ such that
$$
\| \cJ'(u_n) (\varphi(\cdot - z_n)) \| \leq \| \cJ '(u_n)\| \| \varphi(\cdot - z_n)\| \leq D \| \cJ'(u_n) \| \|\varphi\|_{H^1(\R^N)} \to 0.
$$
Therefore
\begin{eqnarray}\label{11}
\int_{\R^N} \nabla v_n \nabla \varphi \, dx + \int_{\R^N} V(x+z_n) v_n \varphi \, dx - \int_{\R^N} g(x,v_n) \varphi \, dx \to 0.
\end{eqnarray}
Since $|z_n| \to \infty$ implies that $V_{loc} (x+z_n) \to 0$ for a.e. $x \in \R^N$, we have 
\begin{eqnarray}\label{22}
\int_{\R^N} (V(x+z_n) - V_{per}(x+z_n))v_n \varphi \, dx = \int_{\R^N} V_{loc} (x+z_n) v_n \varphi \, dx \to 0.
\end{eqnarray}
Combining (\ref{11}) and (\ref{22}) we have that
\begin{eqnarray*}
\cJ_{per}' (v_n)(\varphi) &=& \int_{\R^N} \nabla v_n \nabla \varphi \, dx + \int_{\R^N} V_{per}(x) v_n \varphi \, dx - \int_{\R^N} g(x, v_n) \varphi \, dx = \\
&=& \left[  \int_{\R^N} \nabla v_n \nabla \varphi \, dx + \int_{\R^N} V(x+z_n) v_n \varphi \, dx - \int_{\R^N} g(x, v_n) \varphi \, dx \right] \\
&-& \left[ \int_{\R^N} V_{loc} (x+z_n) v_n \varphi \, dx \right] \to 0,
\end{eqnarray*}
which completes the proof of Step 3. \\
\textbf{Step 4:} \textit{Suppose that there exist $m \geq 1$, $(y_n^k) \subset \mathbb{Z}^N$, $w^k \in H^1 (\R^N)$ for $1\leq k\leq m$ such that
\begin{eqnarray*}
|y_n^k| \to \infty, \ |y_n^k - y_n^{k'}| \to \infty \quad for \ k \neq k', \\
u_n (\cdot + y_n^k) \to w^k \neq 0, \quad for \ each \ 1\leq k \leq m, \\
\cJ_{per}'(w^k)=0, \quad for \ each \ 1\leq k \leq m.
\end{eqnarray*}
Then,
\begin{itemize}
\item[(1)] if $\sup_{z \in \R^N} \int_{B(z,1)} \left| u_n - u_0 - \sum_{k=1}^m w^k (\cdot - y_n^k) \right|^2 \, dx \to 0$ as $n\to\infty$, then
$$
\left\| u_n - u_0 - \sum_{k=1}^m w^k (\cdot - y_n^k) \right\| \to 0;
$$
\item[(2)] if there is $(z_n) \subset \mathbb{Z}^N$ such that 
$$
\liminf_{n \to \infty} \int_{B(z_n, 1+\sqrt{N})} \left| u_n - u_0 - \sum_{k=1}^m w^k (\cdot - y_n^k) \right|^2 \, dx  > 0,
$$
then there is $w^{m+1} \in H^1 (\R^N)$ such that (up to subsequences):
\begin{itemize}
\item[(i)] $|z_n| \to \infty$, \ $|z_n - y_n^k| \to \infty$, \ for \ $1 \leq k \leq m$,
\item[(ii)] $u_n(\cdot + z_n) \rightharpoonup w^{m+1} \neq 0,$
\item[(iii)] $\cJ_{per}'(w^{m+1}) = 0.$
\end{itemize}
\end{itemize}
}
\noindent Suppose that  $\sup_{z \in \R^N} \int_{B(z,1)} \left| u_n - u_0 - \sum_{k=1}^m w^k (\cdot - y_n^k) \right|^2 \, dx \to 0$. Let
$$
\xi_n = u_n - u_0 - \sum_{k=1}^m w^k (\cdot - y_n^k)
$$
and observe that, in view of Lions' lemma, we have $\xi_n \to 0$ in $L^r (\R^N)$. Let us compute
\begin{eqnarray*}
\mathcal{J}'(u_n) (\xi_n) &=& \int_{\mathbb{R}^N} \nabla \xi_n \nabla \xi_n \, dx + \int_{\mathbb{R}^N} \nabla u_0 \nabla \xi_n \, dx + \int_{\mathbb{R}^N} \nabla \left( \sum_{k=1}^m w^k (\cdot - y_n^k) \right) \nabla \xi_n \, dx \\
&+& \int_{\mathbb{R}^N} V(x) \xi_n^2 \, dx + \int_{\mathbb{R}^N} V(x) u_0 \xi_n \, dx \\
&+& \int_{\mathbb{R}^N} V(x) \left( \sum_{k=1}^m w^k (\cdot - y_n^k) \right) \xi_n \, dx - \int_{\mathbb{R}^N} g(x,u_n) \xi_n \, dx.
\end{eqnarray*}
Thus
\begin{eqnarray*}
\| \xi_n \|^2 &=& \mathcal{J}' (u_n) (\xi_n) - \int_{\mathbb{R}^N} \nabla u_0 \nabla \xi_n \, dx - \int_{\mathbb{R}^N} V(x) u_0 \xi_n \, dx - \int_{\mathbb{R}^N} \nabla \left( \sum_{k=1}^m w^k (\cdot - y_n^k )\right) \nabla \xi_n \, dx \\ &-& \int_{\mathbb{R}^N} V(x) \left( \sum_{k=1}^m w^k (\cdot - y_n^k ) \right) \xi_n \, dx + 
\int_{\mathbb{R}^N} g(x,u_n)\xi_n \, dx.
\end{eqnarray*}
Since $\cJ'(u_0)(\xi_n) = 0$, we get
\begin{eqnarray*}
\|\xi_n\|^2 &=& \mathcal{J}'(u_n) (\xi_n) - \int_{\mathbb{R}^N} g(x,u_0) \xi_n \, dx - \sum_{k=1}^m \int_{\mathbb{R}^N} \nabla ( w^k(\cdot - y_n^k)) \nabla \xi_n \, dx \\
&-& \sum_{k=1}^m \int_{\mathbb{R}^N} V_{per} (x) w^k(\cdot - y_n^k) \xi_n \, dx - \sum_{k=1}^m \int_{\mathbb{R}^N} V_{loc} (x) w^k (\cdot - y_n^k) \xi_n \, dx \\
&+& \int_{\mathbb{R}^N} g(x,u_n) \xi_n \, dx.
\end{eqnarray*}
Moreover $\cJ_{per}'(w^k) = 0$, thus
\begin{eqnarray}\label{eq:ksi_n}
\| \xi_n \|^2 &=& \mathcal{J}'(u_n) (\xi_n) - \sum_{k=1}^m \int_{\mathbb{R}^N} g(x, w^k) \xi_n (\cdot + y_n^k) \, dx \\
&-& \sum_{k=1}^m \int_{\mathbb{R}^N} V_{loc}(x) w^k (\cdot - y_n^k) \xi_n \, dx + \int_{\mathbb{R}^N} (g(x,u_n) - g(x,u_0)) \xi_n \, dx.\nonumber
\end{eqnarray}
Observe that
$$
\| \mathcal{J}'(u_n) (\xi_n)\| \leq \| \mathcal{J}'(u_n) \| \|\xi_n\| \to 0.
$$
Since $|\xi_n|_r \to 0$, then by \eqref{eq:NonIneq} and in view of  the H\"older inequality
\begin{eqnarray*}
\left|\int_{\mathbb{R}^N} g(x,u_0) \xi_n \, dx \right| &\leq & \varepsilon \int_{\mathbb{R}^N} |u_0 \xi_n| \, dx + C_\varepsilon \int_{\mathbb{R}^N} |u_0|^{r-1} |\xi_n| \, dx\\
&\to & \varepsilon \limsup_{n\to\infty}\int_{\mathbb{R}^N} |u_0 \xi_n| \, dx
\end{eqnarray*}
 as $n\to\infty$, and taking $\varepsilon\to 0$ we get 
$$\left|\int_{\mathbb{R}^N} g(x,u_0) \xi_n \, dx \right|\to 0.$$
Similarly we show that all integrals in \eqref{eq:ksi_n} tend to $0$.
Thus $\|\xi_n\| \to 0$ and the proof in this case is completed. Suppose that    
$$
\liminf_{n \to \infty} \int_{B(z_n, 1+\sqrt{N})} \left| u_n - u_0 - \sum_{k=1}^m w^k (\cdot - y_n^k) \right|^2 \, dx  > 0,
$$
for some $(z_n) \subset \mathbb{Z}^N$. Then $(i)$ and $(ii)$ hold similarly as in Step 3. To prove $(iii)$, let $v_n = u_n(\cdot + z_n)$. Then for $\varphi \in \cC_0^\infty (\R^N)$ we have
$$
\cJ_{per}'(v_n)(\varphi) - \cJ_{per}'(w^{m+1})(\varphi) \to 0
$$
and  $\cJ_{per}'(v_n)(\varphi) \to 0$, which completes the proof of $(iii)$. \\
\textbf{Step 5:} \textit{Conclusion.} \\
On view of Step 1, we know that $u_n \rightharpoonup u_0$ and $\cJ'(u_0) = 0$, which completes the proof of (a). If condition \eqref{eq:step2} from Step 2 holds, then $u_n \to u_0$ and theorem is true for $\ell = 0$. On the other hand, one has
$$
\liminf_{n\to\infty} \int_{B(y_n,1)} |v_n^1|^2 \, dx > 0
$$
for some $(y_n) \subset \R^N$. For each $y_n \in \R^N$ we will find $z_n \in \mathbb{Z}^N$ such that
$$
B(y_n, 1) \subset B(z_n, 1+\sqrt{N}).
$$
Then
$$
\liminf_{n\to\infty} \int_{B(z_n, 1+\sqrt{N})} |v_n^1|^2 \, dx \geq \liminf_{n\to\infty}\int_{B(y_n,1)} |v_n^1|^2 \, dx > 0.
$$
Therefore in view of Step 3 we find $w$ such that (i)-(iii) hold. Let $y_n^1 = z_n$ and $w^1 = w$. If (1) from Step 4 holds with $m=1$, then (b)-(d) are true. Otherwise (2) holds and we put $(y_n^2) = (z_n)$ and $w^2=w$. Then we iterate the Step 4. To complete the proof of (b)-(d) it is sufficient to show that this procedure will finish after a finite number of steps. Indeed, observe that
$$
\lim_{n\to\infty} \|u_n\|^2-\|u_0\|^2-\sum_{k=1}^m\|w^k\|^2 = \lim_{n\to\infty} \left\| u_n - u_0 - \sum_{k=1}^m w^k (\cdot - y_n^k) \right\|^2 \geq 0
$$
for each $m \geq 1$. Since $w^k$ are critical points of $\cJ_{per}$, there is $\rho_0 > 0$ such that $\|w^k\| \geq \rho_0 > 0$, so after a finite number of steps, say $\ell$ steps, condition (1) in Step 4 will hold. \\
\textbf{Step 6:} \textit{We will show that $(e)$ holds:}
$$
\cJ(u_n)\to \cJ(u_0) + \sum_{k=1}^\ell \cJ_{per} (w^k).
$$
Observe that
\begin{eqnarray*}
&& \mathcal{J}(u_n) = \frac{1}{2} \int_{\mathbb{R}^N} |\nabla u_0 |^2 \, dx + \frac{1}{2} \int_{\mathbb{R}^N} |\nabla (u_n - u_0)|^2 \, dx + \int_{\mathbb{R}^N} \nabla u_0 \nabla (u_n - u_0) \, dx \\
&&+ \frac{1}{2} \int_{\mathbb{R}^N} V (x) u_0^2 \, dx + \frac{1}{2} \int_{\mathbb{R}^N} V(x) (u_n - u_0)^2 \, dx + \int_{\mathbb{R}^N} V(x) u_0 (u_n - u_0) \, dx - \int_{\mathbb{R}^N} G(x,u_n) \, dx
\end{eqnarray*}
and therefore
\begin{eqnarray*}
\mathcal{J}(u_n) &=& \mathcal{J}(u_0) + \mathcal{J}_{per} (u_n - u_0) + \int_{\mathbb{R}^N} \nabla u_0 \nabla (u_n - u_0) \, dx \\
&&+ \frac{1}{2} \int_{\mathbb{R}^N} V_{loc} (x) (u_n - u_0)^2 \, dx + \int_{\mathbb{R}^N} V(x) u_0 (u_n - u_0) \, dx \\
&&+ \int_{\mathbb{R}^N} G(x, u_n - u_0) \, dx + \int_{\mathbb{R}^N} G(x,u_0) \, dx - \int_{\mathbb{R}^N} G(x,u_n) \, dx.
\end{eqnarray*}
It is sufficient to show that
\begin{eqnarray}\label{eq:Step6_1}
&&\int_{\mathbb{R}^N} [ G(x,u_n-u_0) + G(x,u_0) - G(x,u_n)] \, dx \to 0,\\
\label{eq:Step6_2}
&& \mathcal{J}_{per} (u_n - u_0) \to \sum_{k=1}^\ell \mathcal{J}_{per} (w^k).
\end{eqnarray}
Let us consider the function $H : \R^N \times [0,1] \rightarrow \R$ given by $H(x,t)=G(x,u_n-tu_0)$. Therefore
$$
G(x,u_n-u_0)-G(x,u_n) = H(x,1)-H(x,0) = \int_0^1 \frac{\partial H}{\partial s} (x,s) \, ds.
$$
Note that
\begin{eqnarray*}
\int_{\mathbb{R}^N} [ G(x,u_n-u_0) + G(x,u_0) - G(x,u_n)] \, dx &=& \int_{\R^N} \left[ \int_0^1 \frac{\partial H}{\partial s} (x,s) \, ds + G(x,u_0) \right] \, dx \\
&=& \int_{\R^N} \int_0^1 \frac{\partial H}{\partial s} (x,s) \, ds \, dx + \int_{\R^N}G(x,u_0)\, dx\\
&=& \int_0^1 \int_{\R^N} -g(x,u_n-su_0)u_0  \, dx \, ds + \int_{\R^N}G(x,u_0)\, dx.
\end{eqnarray*}
Let $E \subset \R^N$ be a measurable set. From the H\"older inequality we have
\begin{eqnarray*}
\int_{E} |g(x, u_n-su_0)u_0| \, dx &\leq& C \int_E |u_n-su_0| |u_0| \, dx + C \int_E |u_n-su_0|^{r-1} |u_0| \, dx \\
&\leq& C |(u_n-su_0)\chi_E|_2^2 |u_0\chi_E|_2^2 + C|(u_n - su_0)\chi_E|_r^{r-1} |u_0 \chi_E|_r.
\end{eqnarray*}
Then $(g(x, u_n-su_0)u_0)$ is uniformly integrable and by the Vitali convergence theorem we get
$$
\int_0^1 \int_{\R^N} -g(x,u_n-su_0)u_0  \, dx \, ds \to \int_0^1 \int_{\R^N} -g(x,u_0-su_0)u_0  \, dx \, ds.
$$
On the other hand we we have
\begin{eqnarray*}
\int_0^1 \int_{\R^N} -g(x,u_0-su_0)u_0  \, dx \, ds &=&  \int_{\R^N} \int_0^1 -g(x,u_0-su_0)u_0  \, ds \, dx \\
&=& \int_{\R^N} \int_0^1 \frac{\partial}{\partial s} \left[ G(x, u_0-su_0) \right] \, ds \, dx \\
&=& \int_{\R^N} G(x,0)-G(x,u_0) \, dx = \int_{\R^N} -G(x,u_0) \, dx.
\end{eqnarray*}
Finally
$$
\int_{\mathbb{R}^N} [ G(x,u_n-u_0) + G(x,u_0) - G(x,u_n)] \, dx \to \int_{\mathbb{R}^N} [ G(x,u_0) - G(x,u_0)] \, dx = 0,
$$
which completes the proof of \eqref{eq:Step6_1}. In order to show \eqref{eq:Step6_2} observe that
$$
\cJ_{per} (u_n-u_0) = \int_{\R^N} |\nabla (u_n-u_0)|^2 + V_{per}(x)(u_n-u_0)^2 \, dx - \int_{\R^N} G(x, u_n-u_0) \, dx.
$$
Now we show that 
$$
\int_{\R^N} G(x, u_n-u_0) \, dx \to \sum_{k=1}^\ell \int_{\R^N} G(x, w^k) \, dx.
$$
Put $a_m^n = u_n-u_0-\sum_{k=1}^m w^k (\cdot - y_n^k)$. Observe that in \eqref{eq:Step6_1} we have already proved that
$$
\int_{\R^N} G(x,a_0^n) + G(x,u_0)-G(x,u_n) \, dx \to 0
$$
as $n\to\infty$.
Taking $a_0^n(\cdot + y_n^1)$ instead of $u_n$ and $w^1$ instead of $u_0$ we get
\begin{equation}\label{eq:Step6_3}
\int_{\R^N} G(x,a_1^n) + G(x,w^1)-G(x,a_0^n) \, dx \to 0.
\end{equation}
Now taking $a_1^n (\cdot + y_n^2)$ instead of $u_n$ and $w^2$ instead of $u_0$ we get
$$
\int_{\R^N} G(x,a_2^n) + G(x,w^2)-G(x,a_1^n) \, dx \to 0.
$$
Taking into account (\ref{eq:Step6_3}) in the above formula we obtain
$$
\int_{\R^N} G(x,a_2^n)\, dx + \int_{\R^N}G(x,w^2)\, dx+\int_{\R^N}G(x,w^1)\, dx-\int_{\R^N}G(x,a_0^1) \, dx \to 0.
$$
Repeating this reasoning, putting $a^n_{\ell -1} (\cdot + y_n^\ell)$ instead of $u_n$ and $w^\ell$ instead of $u_0$ we get
$$
\int_{\R^N} G(x,a_\ell^n)\,dx + \int_{\R^N} G(x,w^\ell)\, dx - \int_{\R^N} G(x,a^n_{\ell-1})\, dx \to 0.
$$
Using, already proved convergences, we get respectively:
\begin{eqnarray*}
&& \int_{\R^N} G(x,a_\ell^n)\,dx + \int_{\R^N} G(x,w^\ell)\, dx - \int_{\R^N} G(x,a^n_{\ell-1})\, dx \to 0, \\
&& \int_{\R^N} G(x,a_\ell^n)\,dx + \int_{\R^N} G(x,w^\ell)\, dx + \int_{\R^N} G(x,w^{\ell-1})\, dx - \int_{\R^N} G(x, a_{\ell-2}^n) \, dx \to 0,
\\
&& \vdots \\
&& \int_{\R^N} G(x,a_\ell^n) \, dx + \sum_{k=1}^\ell \int_{\R^N} G(x, w^k) \, dx - \int_{\R^N} G(x, a_0^n) \to 0.
\end{eqnarray*}
Observe that $a_\ell^n \to 0$, and therefore
$$
\int_{\R^N} G(x,a_\ell^n) \, dx \to 0.
$$
Finally
$$
\int_{\R^N} G(x, u_n - u_0) \, dx \to \sum_{k=1}^\ell \int_{\R^N} G(x,w^k) \, dx.
$$
Observe that 
$$
\int_{\R^N} \left| \nabla \left( u_n - u_0 - \sum_{k=1}^\ell w^k (\cdot - y_n^k) \right) \right|^2 + V(x) \left( u_n - u_0 - \sum_{k=1}^\ell w^k (\cdot - y_n^k) \right)^2 \, dx \to 0,
$$
which is equivalent to
\begin{eqnarray*}
& & \int_{\R^N} |\nabla (u_n-u_0)|^2 \, dx + \sum_{k=1}^\ell \int_{\R^N} | \nabla w^k (\cdot - y_n^k)|^2 \, dx \\
&& - 2 \int_{\R^N} \sum_{k=1}^\ell \nabla (u_n - u_0) \cdot \nabla w^k (\cdot - y_n^k) \, dx + \sum_{k \neq k'} \int_{\R^N} \nabla w^k (\cdot - y_n^k) \nabla w^{k'}(\cdot - y_n^{k'}) \, dx \\
&& + \int_{\R^N} V_{loc}(x) \left( u_n - u_0 - \sum_{k=1}^\ell w^k (\cdot - y_n^k) \right)^2 \, dx \\
&& + \int_{\R^N} V_{per}(x) \left( u_n - u_0 - \sum_{k=1}^\ell w^k (\cdot - y_n^k) \right)^2 \, dx \to 0.
\end{eqnarray*}
Note that
$$
\int_{\R^N} V_{loc}(x) \left( u_n - u_0 - \sum_{k=1}^\ell w^k (\cdot - y_n^k) \right)^2 \, dx \to 0,
$$
since
$$
\left| \int_{\R^N} V_{loc}(x) \left( u_n - u_0 - \sum_{k=1}^\ell w^k (\cdot - y_n^k) \right)^2 \, dx \right| \leq |V_{loc}|_\infty \left|u_n - u_0 - \sum_{k=1}^\ell w^k (\cdot - y_n^k)\right|_2^2.
$$
Moreover
$$
- 2 \int_{\R^N} \sum_{k=1}^\ell \nabla (u_n - u_0) \cdot \nabla w^k (\cdot - y_n^k) \, dx \to - 2 \sum_{k=1}^\ell \int_{\R^N}|\nabla w^k (\cdot - y_n^k)|^2 \, dx + o(1),
$$
since
$$
\int_{\R^N} \sum_{k=1}^\ell \nabla (u_n - u_0) \cdot \nabla w^k (\cdot - y_n^k) \, dx = \int_{\R^N} \sum_{k=1}^\ell \nabla (u_n (\cdot + y_n^k) - u_0(\cdot + y_n^k)) \cdot \nabla w^k \, dx
$$
and $u_n (\cdot + y_n^k) \rightharpoonup w^k$. Therefore
\begin{eqnarray*}
& & \int_{\R^N} |\nabla (u_n-u_0)|^2 \, dx - \sum_{k=1}^\ell \int_{\R^N} | \nabla w^k (\cdot - y_n^k)|^2 \, dx + \sum_{k \neq k'} \int_{\R^N} \nabla w^k (\cdot - y_n^k) \nabla w^{k'}(\cdot - y_n^{k'}) \, dx \\
&& + \int_{\R^N} V_{per}(x) \left( u_n - u_0 - \sum_{k=1}^\ell w^k (\cdot - y_n^k) \right)^2 \, dx \to 0.
\end{eqnarray*}
Observe that
$$
\int_{\R^N} \nabla w^k (\cdot - y_n^k) \nabla w^{k'}(\cdot - y_n^{k'}) \, dx = \int_{\R^N} \nabla w^k \nabla w^{k'}  (\cdot + y_n^k - y_n^{k'}) \, dx \to 0,
$$
since $|y_n^k - y_n^{k'}| \to \infty$ for $k \neq k'$.  Thus
\begin{eqnarray*}
& & \int_{\R^N} |\nabla (u_n-u_0)|^2 \, dx - \sum_{k=1}^\ell \int_{\R^N} | \nabla w^k (\cdot - y_n^k)|^2 \, dx \\
&& + \int_{\R^N} V_{per}(x) \left( u_n - u_0 - \sum_{k=1}^\ell w^k (\cdot - y_n^k) \right)^2 \, dx \to 0,
\end{eqnarray*}
which is equivalent to
\begin{eqnarray*}
& & \int_{\R^N} |\nabla (u_n-u_0)|^2 \, dx - \sum_{k=1}^\ell \int_{\R^N} | \nabla w^k (\cdot - y_n^k)|^2 \, dx + \int_{\R^N} V_{per}(x) (u_n - u_0)^2 \, dx \\
&& - 2 \sum_{k=1}^\ell \int_{\R^N} V_{per} (x) (u_n - u_0) w^k (\cdot - y_n^k) \, dx + \int_{\R^N} V_{per} (x) \left( \sum_{k=1}^\ell w^k (\cdot - y_n^k) \right)^2 \, dx \to 0.
\end{eqnarray*}
In a similar way we get
$$
- 2 \sum_{k=1}^\ell \int_{\R^N} V_{per} (x) (u_n - u_0) w^k (\cdot - y_n^k) \, dx = -2 \int_{\R^N} V_{per}(x) (u_n - u_0)^2 \, dx + o(1).
$$
Therefore
\begin{eqnarray*}
& & \int_{\R^N} |\nabla (u_n-u_0)|^2 \, dx - \sum_{k=1}^\ell \int_{\R^N} | \nabla w^k (\cdot - y_n^k)|^2 \, dx - \int_{\R^N} V_{per}(x) (u_n - u_0)^2 \, dx \\
&& + \int_{\R^N} V_{per} (x) \left( \sum_{k=1}^\ell w^k (\cdot - y_n^k) \right)^2 \, dx \to 0
\end{eqnarray*}
and
\begin{eqnarray*}
& & \int_{\R^N} |\nabla (u_n-u_0)|^2 \, dx - \sum_{k=1}^\ell \int_{\R^N} | \nabla w^k (\cdot - y_n^k)|^2 \, dx - \int_{\R^N} V_{per}(x) (u_n - u_0)^2 \, dx \\
&& + \sum_{k=1}^\ell \int_{\R^N} V_{per}(x) (w^k (\cdot - y_n^k) )^2 \, dx + \sum_{k \neq k'} \int_{\R^N} V_{per}(x) w^k (\cdot - y_n^k) w^{k'} (\cdot - y_n^{k'}) \, dx \to 0.
\end{eqnarray*}
Note that
$$
\int_{\R^N} V_{per}(x) w^k (\cdot - y_n^k) w^{k'} (\cdot - y_n^{k'}) \, dx \to 0,
$$
and therefore
$$
\cJ_{per} (u_n - u_0) - \sum_{k=1}^\ell \cJ_{per} (w^k) = \cJ_{per} (u_n - u_0) - \sum_{k=1}^\ell \cJ_{per} (w^k(\cdot - y_n^k)) \to 0.
$$
Thus
$$
\cJ_{per} (u_n - u_0) \to \sum_{k=1}^\ell \cJ_{per} (w^k),
$$
which completes the proof of \eqref{eq:Step6_2} and the proof of Theorem \ref{ThDecomposition}.
\end{altproof}

{\bf Acknowledgements.}
The second author was partially supported by the National Science Centre, Poland (Grant No. 2014/15/D/ST1/03638).



\begin{thebibliography}{99}
\baselineskip 2 mm





\bibitem{AlamaLi} S. Alama, Y. Y. Li:
{\em On ''multibump'' bound states for certain semilinear elliptic equations},  Indiana Univ. Math. J. {\bf 41}, (1992), no. 4, 983--1026. 



\bibitem{AmbrosettiRabinowitz} A. Ambrosetti, P. Rabinowitz: {\em Dual variational methods in critical point theory and applications}, J. Functional Analysis {\bf 14} (1973), 349--381.








\bibitem{BartschMederski1} T. Bartsch, J. Mederski:
{\em Ground and bound state solutions of semilinear time-harmonic Maxwell equations in a bounded domain}, Arch. Rational Mech. Anal. {\bf 215} (1), (2015), 283--306.

\bibitem{BartschMederski2} T. Bartsch, J. Mederski:
{\em Nonlinear time-harmonic Maxwell equations in an anisotropic bounded medium}, submitted arXiv:1509.01994


\bibitem{BelmonteBetiaPelinovsky} J. Belmonte-Beitia, D. Pelinovsky: {\em Bifurcation of gap solitons in periodic potentials with a periodic sign-varying nonlinearity coefficient}, Appl. Anal. {\bf 89}, (2010), no. 9, 1335--1350.



\bibitem{BenciGrisantiMeicheletti} V. Benci, C.R. Grisanti, A.M. Micheletti: {\em Existence and non existence of the ground state solution for the nonlinear Schr\"odinger equations with $V(\infty) = 0$}, Topol. Methods in Nonlinear Anal. {\bf 26}, (2005), 203--219.



\bibitem{NonKerrBook} A. Biswas, S. Konar: {\em Introduction to non-Kerr Law Optical Solitons}, Chapman and Hall (2006).



\bibitem{Buryak} A.V. Buryak, P. Di Trapani, D.V. Skryabin, S. Trillo: {\em Optical solitons due to quadratic nonlinearities: from basic physic to futuristic applications}, Physics Reports {\bf 370}, (2002) 63--235.


\bibitem{CostaTehraniCalVPDE} D. Costa, H. Tehrani: {\em Existence of positive solutions for a class of indefinite elliptic problems in $\R^N$}, Cal. Var. {\bf 13} (2), 159--189. 


\bibitem{CostaTehrani} D. Costa, H. Tehrani: {\em Existence and multiplicity results for a class of Schr\"{o}dinger equations with indefinite nonlinearities}, Adv. Differential Equations {\bf 8}, (2003) 1319--1340.


\bibitem{CotiZelati} V. Coti-Zelati, P. Rabinowitz: {\em Homoclinic type solutions for a semilinear elliptic PDE on $\R^n$}, Comm. Pure Appl. Math. {\bf 45}, (1992), no. 10, 1217--1269.



\bibitem{Doerfler} W. D\"orfler, A. Lechleiter, M. Plum, G. Schneider, C. Wieners, {\em Photonic Crystals: Mathematical Analysis and Numerical Approximation}, Springer Basel (2012)


\bibitem{Figueiredo} G. Figueiredo, H. R. Quoirin: {\em Ground states of elliptic problems involving non homogeneous operators}, to appear in Indiana Univ. Math. J. {\bf 65} (3), (2016), arXiv:1503.07479.


\bibitem{GuoMederski} Q. Guo, J. Mederski: {\em Ground states of nonlinear Schr\"odinger equations with sum of periodic and inverse-square potentials},
J. Differential Equations {\bf 260}, (2016) 4180--4202.

\bibitem{JeanjeanTanakaIndiana2005} L. Jeanjean, K. Tanaka: {\em A positive solution for a nonlinear Schr\"odinger equation on $\R^N$}, Indiana Univ. Math. Journal {\bf 54} (2), (2005), 443--464.

\bibitem{Kuchment} P. Kuchment: {\em The mathematics of photonic crystals}, Mathematical modeling in optical science, Frontiers Appl. Math., {\bf 22}, SIAM, Philadelphia (2001), 207--272.



\bibitem{LiWangZeng} Y. Li, Z.-Q. Wang, J. Zeng: {\em Ground states of nonlinear Schr\"odinger equations with potentials}, Ann. Inst. H. Poincaré Anal. Non Linéaire 23 (2006), no. 6, 829--837. 



\bibitem{Lions84} P.L. Lions: {\em The concentration-compactness principle in the calculus of variations. The locally compact case. Part I and II}, Ann. Inst. H. Poincar\'e, Anal. Non Lin\'e are., {\bf 1}, (1984), 109--145; and
223--283.

\bibitem{LiuYang} F. Liu, J. Yang: {\em Nontrivial solutions of Schr\"{o}dinger equations with indefinite nonlinearities}, J. Math. Anal. Appl. {\bf 334}, (2007), 627--645.

\bibitem{Liu} S. Liu: {\em On superlinear Schr\"odinger equations with periodic potential}, Calc. Var. Partial Differential Equations {\bf 45} (2012), no. 1-2, 1--9.



\bibitem{MederskiTMNA2014} J. Mederski: {\em Solutions to a nonlinear Schr\"odinger equation with periodic potential and zero on the boundary of the spectrum}, Topol. Methods Nonlinear Anal. {\bf 46} (2), (2015), 755--771.

\bibitem{MederskiNLS2014} J. Mederski: {\em Ground states of a system of nonlinear Schr\"odinger equations with periodic potentials},  to appear in Comm. Partial Differential Equations (2016), DOI:10.1080/03605302.2016.1209520.




\bibitem{Pankov} A. Pankov: {\em Periodic Nonlinear Schr\"odinger Equation with Application to Photonic Crystals}, Milan J. Math. {\bf 73}, (2005), 259--287.

\bibitem{PankovDecay} A. Pankov: {\em On decay of solutions to nonlinear Schr\"odinger equations},  Proc. Amer. Math. Soc. {\bf 136}, (2008), 2565--2570.

\bibitem{PankovNotes} A. Pankov: {\em Lecture Notes on Schr\"odinger Equations}, Nova Publ., 2007




\bibitem{Rabinowitz:1992} P.H. Rabinowitz: {\em On a class of nonlinear Schr\"odinger equations}, Z. Angew. Math. Phys. {\bf 43}, (1992), 270--291.

\bibitem{ReedSimon} M. Reed, B. Simon: {\em Methods of Modern Mathematical Physics, Analysis of Operators, Vol. IV}, Academic Press, New York, 1978.


\bibitem{Simon} B. Simon: {\em Schr\"odinger semigroups}, Bull. Amer. Math. Soc. {\bf 7} (3), (1982), 447--526.



\bibitem{NonlinearPhotonicCrystals} R. E. Slusher, B. J. Eggleton: {\em Nonlinear Photonic Crystals}, Springer 2003.

\bibitem{Struwe} M. Struwe: {\em Variational Methods}, Springer 2008.



\bibitem{SzulkinWeth} A. Szulkin, T. Weth: {\em Ground state solutions for some indefinite variational problems}, J. Funct. Anal. {\bf 257}, (2009), no. 12, 3802--3822. 

\bibitem{SzulkinWethHandbook} A. Szulkin, T. Weth: {\em The method of Nehari manifold. Handbook of nonconvex analysis and applications}, Handbook of nonconvex analysis and applications, 597--632, Int. Press, Somerville, 2010.




\bibitem{Willem} M. Willem: {\em Minimax Theorems}, Birkh\"auser Verlag 1996.



\end{thebibliography}
\end{document}